\documentstyle[dina4]{article}
\newcommand{\N}{{\Bbb N}}

\newcommand{\C}{{\Bbb C}}

\renewcommand{\P}{{\Bbb P}}

\newcommand{\F}{{\Bbb F}}

\newcommand{\ka}{{\cal A}}

\newcommand{\kc}{{\cal C}}
\newcommand{\kd}{{\cal D}}

\newcommand{\kf}{{\cal F}}
\newcommand{\kg}{{\cal G}}
\newcommand{\kh}{{\cal H}}

\newcommand{\km}{{\cal M}}
\newcommand{\kn}{{\cal N}}
\newcommand{\ko}{{\cal O}}

\newcommand{\ks}{{\cal S}}
\newcommand{\kt}{{\cal T}}
\newcommand{\ku}{{\cal U}}

\input amssym.def
\input amssym
\begin{document}
\newtheorem{lemma}{Lemma}[section]
\newtheorem{remark}[lemma]{Remark}
\newtheorem{proposition}[lemma]{Proposition}
\newtheorem{example}[lemma]{Example}
\newtheorem{theorem}[lemma]{Theorem}
\newtheorem{corollary}[lemma]{Corollary}
\newtheorem{definition}[lemma]{Definition}
\newtheorem{proof}[lemma]{Proof}
\newtheorem{sub}[lemma]{}

\title{{\bf Equianalytic and equisingular families of curves on surfaces}}
\author{
Gert-Martin Greuel\\
Universit\"at Kaiserslautern\\
Fachbereich Mathematik\\
Erwin-Schr\"odinger-Stra\ss e\\
D -- 67663 Kaiserslautern
\and
Christoph Lossen\\
Universit\"at Kaiserslautern\\
Fachbereich Mathematik\\
Erwin-Schr\"odinger-Stra\ss e\\
D -- 67663 Kaiserslautern}

\maketitle
\thispagestyle{empty}
\date{}
\vspace{2cm}

\tableofcontents
\vspace{1cm}
\setcounter{page}{0}
\newpage

\section*{Introduction}\addcontentsline{toc}{section}{Introduction}

We consider flat families of reduced curves on a smooth surface $S$ such that
for each member $C$ of the family the number of singular points of $C$ and for
each singular point $x \in C$ the ``singularity type'' of $(C,x)$ is fixed.
Fixing these data imposes conditions on the space of all curves and we obtain
in this way a locally closed subscheme of the Hilbert  scheme $H_S$ of $S$.
We are mainly concerned with the study of the  equianalytic
$(H^{ea}_S)$ respectively the equisingular  Hilbert scheme $(H^{es}_S)$, which
are defined
by fixing the analytic respectively the (embedded) topological type of the
singularities.   We show that fixing the analytical (respectively topological)
type of $(C,x)$ imposes, at most, $\tau (C,x) =$ Tjurina number of $(C,x)$
(respectively, at most, $ \mu (C,x)-$ mod$(C,x)$, where $\mu$ denote the
Milnor number, mod  the
 modality in the sense of \cite{AGV}) conditions with equality if
$H^1(C,\kn^{ea}_{C/S})$ (respectively $H^1(C,\kn^{es}_{C/S})$) vanish.   Here
$\kn^{ea}_{C/S}$ (respectively $\kn^{es}_{C/S}$) denote the equianalytic
(respectively equisingular) normal bundle.   The vanishing of $H^1$ implies
the independence of the imposed conditions and the smoothness of $H^{ea}_S$
(respectively $H^{es}_S$) at $C$ (cf.\ \S 3).

In Theorem 3.7 we prove sufficient conditions for the vanishing of
$H^1(C,\kn^{ea}_{C/S})$ (respectively $H^1(C,\kn^{es}_{C/S})$); for the
special case $S = \P^2$ we obtain an additional criterion in Corollary 3.12.
For  the proof we use a
vanishing theorem of \cite{GrK} which is an improvement upon the usual
vanishing theorem for sheaves which are not locally free.   The local
isomorphism defect isod$_x (\kn^{ea}_{C/S}, \ko_C)$, which is introduced in
3.4, measures how much $\kn^{ea}_{C/S}$ fails to be free at $x$, and similar
for $\kn^{es}_{C/S}$.   In many cases of interest, in particular for
$\kn^{es}_{C/S}$ and related sheaves, the isomorphism defect is quite big and
gives a considerable improvement of the desired vanishing results.
Therefore, we make some effort to compute respectively estimate it for certain
classes of singularities in \S 4.   In \S 5 we give some explicit examples and
applications.

The present work is, in some sense, a continuation of some part of \cite{GrK},
where only equianalytic families were considered.
   Our results about the smoothness of $H^{es}_{\P^2}$, which are valid for
arbitrary singularities,  contain the previously known
facts for curves with only ordinary multiple points (cf.\ \cite{Gia}) as a
special case;  concerning the smoothness of $H^{ea}_{\P^2}$ they are an
improvement of \cite{Sh1}.   For concrete applications of the theorems of this
paper it is important to have good estimates for the isomorphism defects.
Most of the formulas concerning these, together with further refinements and
detailed
proofs, appeared in \cite{Lo}.   Our results for $\P^2$ seem to be quite sharp
for small $d$ but are asymptotically weaker than those of Shustin \cite{Sh3}
which are quadratic in $d$.   However, the methods presented here work for
arbitrary surfaces and may be
combined with Shustin's to provide asymptotically optimal results for curves
on (some classes of) rational surfaces.   This will be the subject of a
forthcoming joint paper.
\newpage

\section{Equisingular deformations of plane curve singularities}

In this paragraph we recall some definitions and results due to J.\ Wahl in
the framework of formal deformation theory (cf.\ \cite{Wa2}), transfer them to
the complex analytic category and obtain some additional results which are
used later.

\begin{sub}\label{1.1}{\rm
Let $(C,0) \subset (\C^2,0)$ be a reduced plane curve singularity and $f \in
\ko_{\C^2,0} = \C\{u,v\}$ a convergent power series defining the germ $(C,0)$.
Furthermore, let $m = \mbox{ mult}_0(C)$ denote the multiplicity of $(C,0)$,
that is $f \in \frak{m}^m_{\C^2,0}\backslash\frak{m}^{m+1}_{\C^2,0}$ where
$\frak{m}_{X,x}$ denotes the maximal ideal of a germ $(X,x)$.   Consider a
deformation $\varphi : (\kc,0) \to (T,0)$ of $(C,0)$ over an arbitrary complex
germ $(T,0)$ together with a section $\sigma : (T,0) \to (\kc,0)$.   Without
loss
of generality we may assume $\varphi$ to be embedded, that is $\varphi$ is
given by a commutative diagram

\unitlength1cm
\begin{picture}(5,2.5)
\put(3.5,2){$(C,0)$}
\put(4,1.75){\vector(0,-1){1}}
\put(3.9,0.25){0}

\put(4.75,2){$\hookrightarrow$}
\put(4.75,0.25){$\in$}

\put(5.5,2){$(\kc,0)$}
\put(5.75,0.75){\vector(0,1){1}}
\put(6,1.75){\vector(0,-1){1}}
\put(5.4,1.25){$\sigma$}
\put(6.15,1.25){$\varphi$}

\put(5.5,0.25){$(T,0)$}

\put(6.75,2){$\hookrightarrow$}
\put(7.5,2){$(\C^2 \times T,0)$}

\put(7.8,1.75){\vector(-3,-2){1.5}}
\put(7,1){$pr$}
\end{picture}

where $pr$ is the (natural) projection, $\sigma$ maps to the trivial section
and $(\kc,0)$ is a hypersurface germ of
$(\C^2 \times T,0)$ defined by a power series $F \in \ko_{\C^2 \times T,0}$.
Let $I_{\sigma(T)}$ denote the ideal of  $\sigma(T,0) \subset (\C^2 \times
T,0)$, then we call the deformation with section $(\varphi,\sigma)$ {\sl
equimultiple}, if $F \in I^m_{\sigma(T)}$ (which is, of course, independent of
the choice of the embedding and the choice of $F$).
}
\end{sub}

\begin{sub}\label{1.2}{\rm
Before defining equisingular deformations, let us recall that the {\sl
equisingularity type}\/ (or {\sl topological type}) of $(C,0)$ may be defined
as
follows:  consider an embedded resolution of $(C,0) \subset (\C^2,0)$ given by
a sequence of blowing up points $(N \ge 1)$:
\begin{equation}
M_N \buildrel\pi_N\over \to M_{N-1} \to \cdots \to M_1 \buildrel\pi_1\over\to
M_0 = (\C^2,0).
\end{equation}

Let $C_i \subset M_i$ be the strict and $\hat{C}_i \subset M_i$ the reduced
total transform of $C_0 := (C,0) \subset M_0$ under $\psi_i := \pi_1 \circ
\cdots \circ \pi_i$.   Assume that

\begin{itemize}
\item $\pi_1$ blows up $0$;
\item for $i = 2, \ldots, N$, $\pi_i$ blows up singular points of $\hat{C}_i$;
\item $\hat{C}_N$ has only singularities of type $A_1$.
\end{itemize}
Hence, $C_N$ is smooth and $\psi_N$ induces a resolution of $(C,0)$.   If we
choose the minimal resolution (that is blowing up only non--nodes of
$\hat{C}_i$) we obtain a well--defined system of multiplicity sequences (cf.\
\cite{BrK}), which defines the equisingularity type of $(C,0)$.   It is
well--known (and was proved by Zariski in \cite{Zar}) that this system of
multiplicity sequences determines the embedded topological type of $(C,0)$ and
vice versa (cf.\ \cite{BrK}, 8.4 and \cite{Zar} for further
characterizations).
}
\end{sub}

\begin{sub}\label{1.3}{\rm
Consider a deformation $\varphi : (\kc,0) \to (T,0)$ of $(C,0)$ (without
section) and assume it to be embedded

\unitlength1cm
\begin{picture}(4,2.5)
\put(4.5,2){$(\kc,0)$}
\put(5,1.75){\vector(0,-1){1}}
\put(4.7,1.25){$\varphi$}
\put(4.5,0.25){$(T,0)$}

\put(5.75,2){$\hookrightarrow$}
\put(6.5,2){$(\C^2 \times T,0)$}

\put(6.75,1.75){\vector(-3,-2){1.5}}
\put(6.2,1.15){$pr$}
\end{picture}

$\varphi$ is called {\sl equisingular}\/ (cf.\ \cite{Wa2}, \S 3, \S 7) if
there exists a sequence of blowing up subspaces
\begin{equation}
\km_N \buildrel\tilde{\pi}_N\over\to \km_{N-1} \to \cdots \to \km_1
\buildrel\tilde{\pi}_1\over \to \km_0 = (\C^2 \times T,0)
\end{equation}
such that if $\kc_i \subset \km_i$ denotes the strict and $\hat{\kc}_i \subset
\km_i$ the reduced total transform of $\kc_0 = (\kc,0)$ under $\tilde{\psi}_i
:= \tilde{\pi}_1 \circ \cdots \circ \tilde{\pi}_i$, the following holds:

\begin{itemize}
\item[(i)] sequence (1) is induced by (2) via the base change $0 \mapsto T$;

\item[(ii)] there is a section $\sigma : (T,0) \to \kc_0$ of $\varphi_0 =
\varphi$ such that $\varphi$ is equimultiple along $\sigma$ and
$\tilde{\pi}_1$ blows up $\sigma(T,0) \subset \km_0$;

\item[(iii)] for $i = 1, \ldots, N$ there are sections $(T,0) \to \kc_i$ of
$\varphi_i = \varphi \circ \tilde{\psi}_i \mid_{\kc_i} : \kc_i \to (T,0)$
through all singular points of $\hat{C}_i$ (each of those sections being
mapped via $\tilde{\pi}_i$ to such a section of $\varphi_{i-1}$) such that
$\varphi_i$ is equimultiple along them.   $\tilde{\pi}_{i+1}: \km_{i+1}
\to \km_i$ blows up the sections going through those singular points of
$\hat{C}_i$ which are
blown up by $\pi_{i+1}$ ($i \le N-1$).
\end{itemize}

The sections of (ii), (iii) are called a (compatible) system of {\sl
equimultiple sections} of $\varphi$ through all infinitely near points of
$(C,0)$.

This definition is obviously independent of the embedding of $\varphi$.
Moreover, since an equimultiple deformation of an $A_1$--singularity is
trivial, the definition is also independent of the embedded resolution (1).
The section $\sigma$ of (ii) is called a {\sl singular section}\/ of
$\varphi$.   If $(T,0)$ is reduced, then $\varphi : (\kc,0) \to (T,0)$ is
equisingular if and only if for a small good representative $\varphi : \kc \to
T$ and for all $t \in T$ there exists an $x \in \varphi^{-1}(t)$ such that the
Milnor number $\mu(\varphi^{-1}(t),x)$ is equal to the Milnor number $\mu
(C,0)$.   The existence of a singular section was shown by B.\ Teissier
(\cite{Te}, \S 5).
}
\end{sub}

The following theorem is basically due to J.\ Wahl (\cite{Wa2}, Theorem 7.4):

\begin{theorem} Let $\varphi : (\kc,0) \to (T,0)$ be any equisingular
deformation of the reduced plane curve singularity $(C,0) \subset (\C^2,0)$.

\begin{itemize}
\item[(i)] The equimultiple sections through all infinitely near points of
$(C,0)$ which
are required to exist for $\varphi$ are uniquely determined.

\item[(ii)] Let $\phi : \kc_{(C,0)} \to S_{(C,0)}$ be the semiuniversal
deformation of $(C,0)$.   Then there exists a smooth subgerm $S^{es}_{(C,0)}
\subset S_{(C,0)}$ such that if $\varphi$ is induced from $\phi$ via the base
change $\psi : (T,0) \to S_{(C,0)}$, then $\psi$ factors through
$S^{es}_{(C,0)}$.   In particular, the restriction of $\phi$ to
$S^{es}_{(C,0)}$ is a {\rm semiuniversal equisingular deformation} of
$(C,0)$.

\item[(iii)] Let $T_\varepsilon := \mbox{ Spec}(\C[\varepsilon]/\varepsilon^2)$
be the
base space of first order infinitesimal deformations.   The set

\begin{center}
$I^{es}(C,0) := \{g \in \C\{u,v\}\; \Big|$ \raisebox{1.5ex}{$F = f +
\varepsilon g$
defines an equisingular}\hspace{-5.25cm}\raisebox{-1.5ex}{deformation of
$(C,0)$
over $T_\varepsilon$}\hspace{1cm}$\Big\}$
\end{center}

is an ideal, the {\rm equisingularity ideal}\/ of $(C,0)$.   Especially it
contains the Jacobian ideal
\[
j(C,0) = (f,\frac{\partial f}{\partial u},\; \frac{\partial f}{\partial v})
\cdot \C \{u,v\}
\]
and the vector space $I^{es}(C,0)/j(C,0)$ is isomorphic to the tangent
space of $S^{es}_{(C,0)} \subset S_{(C,0)}$.
\end{itemize}
\end{theorem}

{\bf Proof}:  Wahl considers only deformations over Artinian spaces $(T,0)$
but the above facts follow easily from his results:

\begin{itemize}
\item[(i)] Since we require the existence of holomorphic sections over
arbitrary complex germs $(T,0)$, by Wahl these are unique modulo arbitrary
powers of the maximal ideal $\frak{m}_{T,0}$, hence unique.

\item[(ii)] The existence of a smooth formal semiuniversal equisingular
deformation of $(C,0)$ was proved by Wahl.   The existence of a convergent
representative can be deduced from his result by applying Artin's and Elkik's
algebraization theorems.   A simple direct proof, using the deformation of the
parametrization, is given in \cite{Gr}.

\item[(iii)] follows directly from (\cite{Wa2}, Proposition 6.1).\hfill $\Box$
\end{itemize}

\begin{proposition}\label{1.5}  Openness of versality holds for equisingular
deformations, that is if $\varphi : (\kc,0) \to (T,0)$ is an equisingular
deformation of $(C,0)$, then for any equisingular representative $\varphi :
\kc \to T$ together with the singular section $\sigma : T \to \kc$ the set of
 points $t \in T$ such that $(\kc, \sigma(t)) \to (T,t)$ is a versal
deformation of $(\varphi^{-1}(t), \sigma(t))$ is a Zariski--open subspace of
$T$.
\end{proposition}

{\bf Proof}:  This follows quite formally from a criterion for openness of
versality due to Flenner (\cite{Fl}, Satz 4.3).\hfill $\Box$

\begin{sub}\label{1.6}{\rm
Let $\mu(c,0) = \dim_\C(\C\{u,v\}/(\frac{\partial f}{\partial u},
\frac{\partial f}{\partial v}))$ respectively $\tau(C,0) =
\dim_\C(\C\{u,v\}/(f, \frac{\partial f}{\partial u},\; \frac{\partial
f}{\partial v}))$ denote the Milnor respectively Tjurina number of $(C,0)$.
It is well--known that a
deformation of $(C,0)$ over a reduced base $(T,0)$ is equisingular if and only
if the Milnor number is constant along the (unique) singular section.   Hence
$S^{es}_{(C,0)}$, being smooth,  coincides with the $\mu$--constant stratum of
$S_{(C,0)}$.
The codimension of $S^{es}_{(C,0)}$ in $S_{(C,0)}$ is (by Theorem 1.4
(ii) and (iii)) equal to
\[
\tau^{es}(C,0) = \dim_\C(\C\{u,v\}/I^{es}(C,0)).
\]

Together with a result of Gabrielov (\cite{Gab}), which states that the {\sl
modality}\/ mod$(f)$ of the function $f$ with respect to right equivalence
(cf.\ \cite{AGV}) is equal to the dimension of the $\mu$--constant stratum of
$f$ in the ($\mu$--dimensional) semiuniversal unfolding of $f$, we obtain the
following
}
\end{sub}

\begin{lemma}
For any reduced plane curve singularity $(C,0)$ defined by $f
\in \C\{u,v\}$, we have
\[
\tau^{es}(C,0) = \mu(C,0) - \mbox{ mod}(f).
\]
\end{lemma}

\begin{sub}\label{1.9}{\rm
If $\sim$ denotes any equivalence relation of plane curve singularities, a
$\sim$--{\sl singularity type} is by definition a (not ordered) tuple ${\cal S}
=
((C_1,x_1)/\sim, \ldots, (C_m,x_m)/\sim)$ of $\sim$-equivalence classes with
$m$ a non--negative integer.   In this paper we are mainly interested in the
following two
cases:

\begin{itemize}
\item $\sim =$ analytic equivalence (isomorphism of complex space germs), in
this
case we call the $\sim$--singularity type {\sl analytic type} and denote it by
$\ka$.

\item $\sim =$ topological equivalence (embedded homeomorphism of complex
space germs
(cf.\ 1.3)), the corresponding singularity type is called {\sl
equisingularity type} or {\sl topological type} and denoted by $\kt$.
\end{itemize}

If $(C,x)$ is a reduced plane curve singularity, then ${\cal S}(C,x) =
(C,x)/\sim$
denotes its singularity type and if $C$ is a reduced curve with finitely many
singular points $x_1, \ldots, x_m$ which are all planar, then ${\cal S}(C) =
((C,x_1)/\sim, \ldots, (C,x_m)/\sim)$ is the $\sim$--singularity type of $C$.
For $\ks = \ka$ we obtain $\ka(C)$, the {\sl equianalytic type} of $C$, and
for $\ks = \kt$ we obtain $\kt(C)$, the {\sl equisingular type} of $C$.

As equisingular deformations preserve the topological type, the equianalytic
deformations preserve the analytic type of each fibre, where a deformation
$\varphi : (\kc,0) \to (T,0)$ of $(C,0)$ is called {\sl equianalytic} if
$(\kc,0)$ is analytic isomorphic to $(C \times T,0)$ over $(T,0)$, that is,
$\varphi$ is analytically trivial.}
\end{sub}
\newpage

\section{The equianalytic and equisingular Hilbert scheme}

\begin{sub}\label{2.1}{\rm
Let $S$ be a smooth surface, $T$ a complex space, then by a {\sl family of
embedded (reduced) curves over} $T$ we mean a commutative diagram
\[
\begin{array}{lcl}
\kc & \buildrel j\over\hookrightarrow & \;\;S \times T\\
\varphi\searrow &  & \swarrow pr\\
& T &
\end{array}
\]
where $\varphi$ is a proper and flat morphism such that for all
points $t \in T$ the fibre $\varphi^{-1}(t)$ is a {\sl curve}\/ (that is a
reduced pure 1--dimensional complex space), moreover, $j : \kc
\hookrightarrow S \times T$ is a closed embedding and $pr$ denotes the natural
projection.   Such a family is called {\sl equianalytic} (respectively {\sl
equisingular}) if for all $t \in T$ the induced (embedded) deformation of each
singular point of $\varphi^{-1}(t)$ over ($T,t$) is equianalytic (respectively
equisingular) --- along the unique singular section $\sigma$ (cf.\ \ref{1.3}).
}
\end{sub}

\begin{sub}\label{2.2}{\rm
The  Hilbert functor $\kh ilb_S$ on the category of complex spaces defined by
\[
\kh ilb_S(T) := \{\mbox{ subspaces of } S \times T, \mbox{ proper and flat over
}
T\}
\]
is well--known to be representable by a complex space $H_S$ (cf.\ \cite{Bin}).
This means there is a universal family
\[
\begin{array}{lcl}
\ku & \buildrel j\over\hookrightarrow & \;\; S \times H_S\\
\varphi\searrow & & \swarrow pr\\
& H_S &
\end{array}
\]
such that each element of $\kh ilb_S(T)$, $T$ a complex space, can be induced
from
$\varphi$ via base change by a {\sl unique} map $T \to H_S$.
We define the {\sl equianalytic} (respectively {\sl equisingular}) {\sl
Hilbert functor} $\kh ilb^{ea}_S$ (respectively ${\kh}ilb^{es}_S$) to be the
subfunctor of $\kh ilb_S$ with

\[
\begin{array}{lcl}
 \kh ilb^{ea}_S(T) & := &  \{\mbox{equianalytic families of embedded curves
over }T\}\\[0.5ex]
 \kh ilb^{es}_S(T) & := & \{\mbox{equisingular families of embedded curves
over } T\}
\end{array}
\]
Moreover, fixing the analytic (respectively topological) singularity type
(cf.\ 1.8), we define
\[
\begin{array}{lcl}
 \kh ilb^\ka_S(T) & := & \{\mbox{families in }  \kh ilb^{ea}_S(T)
\mbox{ whose  fibres have (the fixed) analytic singularity type }
\ka\}\\[0.5ex]
 \kh ilb^\kt_S(T) & := & \{ \mbox{families in }  \kh ilb^{es}_S(T) \mbox{ whose
 fibres have (the fixed) topological singularity type } \kt\}
\end{array}
\]
}
\end{sub}

\begin{proposition}
Let $\ka$ be an analytic singularity type corresponding to the topological
type $\kt$, then the functors $\kh ilb^\ka_S$ and $\kh ilb^\kt_S$ are
representable
by locally closed subspaces $H^\ka_S \subset H^\kt_S \subset H_S$.
\end{proposition}

\begin{corollary}
The functor $\kh ilb^{ea}_S$ (respectively $\kh ilb^{es}_S$) is representable
by a
complex space $H^{ea}_S$ (respectively $H^{es}_S$) which is given as the
disjoint union of all $H^{\ka}_S$ (respectively $H^\kt_S$).
\end{corollary}

\begin{remark}{\rm  If $S = \P^2$ and if we fix the degree of all fibres to be
$d$, then
$H^{es,d}_{\P^2} \subset \P^N$ (where $N = \frac{d^2 + 3d}{2}$) is given as a
{\sl finite}\/ disjoint union of locally closed subspaces, while in general
$H^{ea,d}_{\P^2}
\subset \P^N$ is  an {\sl infinite}\/ union.
}
\end{remark}

\begin{sub}\label{2.3}{\rm
For the proof of Proposition 2.3 we need the following Lemma, which, for
the equianalytic case, is proven in (\cite{GrK}, Lemma 1.4).   A proof for
the equisingular case  is given in  \cite{Gr}.
}
\end{sub}

\begin{lemma}
Let $(C,x)$ be the germ of an isolated plane curve singularity and $\varphi
:(\kc,x) \to (B,b)$ a deformation of $(C,x)$, then there are unique closed
subgerms $(B^{ea},b) \subset (B^{es},b) \subset (B,b)$ such that for any
morphism $f : (T,t) \to (B,b)$:
\[
\begin{array}{l}
f^\ast\varphi \mbox{ is an equianalytic deformation if and only if } f(T,t)
\subset (B^{ea},b)\\
f^\ast\varphi \mbox{ is an equisingular deformation if and only if } f(T,t)
\subset (B^{es},b)
\end{array}
\]

Moreover, if $\phi :\kc_{(C,x)} \to S_{(C,x)}$ denotes the semiuniversal
deformation of $(C,x)$ and if $\psi : (B,b) \to S_{(C,x)}$ is any morphism
inducing $\varphi$ via pull--back, then $(B^{ea},b) = (\psi^{-1}(0),b)$ and
$(B^{es},b) = (\psi^{-1}(S^{es}_{(C,x)}),b)$.
\end{lemma}

\begin{sub}\label{2.4}{\rm
{\bf Proof of Proposition 2.3}:  First we have to remark that the
condition for all fibres to be reduced curves defines an open subspace
$\widetilde{H}_S \subset H_S$.   Now, let $b \in \widetilde{H}_S$ be such that
the fibre $\varphi^{-1}(b)$ in the universal family has topological type
$\kt$, then by Lemma 2.7 for each $x \in \varphi^{-1}(b)$ there is a
unique closed subspace $(H_x^{es},b) \subset (\widetilde{H}_S,b)$ such that a
morphism $f : (T,t) \to (\widetilde{H}_S,b)$ factors through $(H^{es}_x,b)$ if
and only if $f^\ast\varphi$ is an equisingular deformation.

Let
\[
(H^{es},b) := \bigcap\limits_{x \in \varphi^{-1}(b)} (H_x^{es},b) \subset
(\widetilde{H}_S,b)
\]
and $H^{es}(b) \subset \widetilde{H}_S$ be a small (unique) representative,
then $\cup H^{es}(b)$, where the union is taken over all $b$ whose fibre
$\varphi^{-1}(b)$ has topological type $\kt$, defines a locally closed subspace
of $H_S$ which obviously represents
$\kh ilb^\kt_S$.   The statement for $\kh ilb^\ka_S$ follows in the
same manner. \hfill$\Box$
}
\end{sub}
\newpage

\section{Completeness of the equianalytic and equisingular characteristic
linear series}

\begin{sub}\label{3.1}{\rm
Let $S$ be a smooth complex surface and $C \subset S$ a reduced compact curve.
Then a {\sl deformation of} $C/S$ over the pointed complex space $T$, $0 \in
T$,  is
a triple $(\kc, \tilde{i}, j)$ such that we obtain a Cartesian diagram

\unitlength1cm
\begin{center}
\begin{picture}(7,2.5)
\put(2.5,2){$C$}
\put(3,2){$\hookrightarrow$}
\put(3.15,2.2){$j$}
\put(3.75,2){$\kc$}
\put(4.3,2.1){\line(1,0){0.5}}

\put(2.3,1.5){$i\;\cap$}
\put(3.75,1.5){$\cap\;\tilde{i}$}

\put(2.5,1){$S$}
\put(3,1){$\hookrightarrow$}
\put(3.75,1){$S \times T$}
\put(5,1){$\Phi$ flat}

\put(2.5,0.5){$\downarrow$}
\put(3.75,0.5){$\downarrow\; \pi$}

\put(2.5,0){$0$}
\put(3.2,0){$\in$}
\put(3.75,0){$T$}
\put(4.8,0.1){\vector(-1,0){0.5}}

\put(4.8,0.1){\line(0,1){2}}

\end{picture}
\end{center}

where $j$ is a closed embedding and the composed morphism $\Phi := \pi \circ
\tilde{i}$ is flat ($S \hookrightarrow S \times T$ denotes the canonical
embedding with image $S \times \{0\}$ and $\pi$ is the projection).   Two
deformations $(\kc, \tilde{i}, j)$ and $(\kc', \tilde{i}', j')$ of $C/S$
over $T$ are {\sl isomorphic} if there exists an isomorphism $\kc \simeq
\kc'$ such that the obvious diagram (with the identity on $S \times T$)
commutes.   $\kd e\!f_{C/S}$ denotes the deformation functor from pointed
complex
spaces to sets defined by
\[
\kd e\!f_{C/S} (T) := \{\mbox{isomorphism classes of deformations of } C/S
\mbox{
over } T\}
\]
and we have the natural forgetful morphism $\kd e\!f_{C/S} \to \kd e\!f_C$
given by
$(\kc, \tilde{i}, j) \mapsto$ \hbox{$(\kc, \Phi = \pi \circ \tilde{i}, j)$},
where
$\kd e\!f_C$ denotes the functor of isomorphism classes of deformations of $C$
(forgetting the embedding).   Furthermore, for each point $x \in C$, we
consider the morphism $\kd e\!f_C \to \kd e\!f_{C,x}$ where $\kd e\!f_{C,x}$
denotes the
functor of isomorphism classes of deformations of the analytic germ
$(C,x)$.

Let $T_\varepsilon =$ Spec $(\C[\varepsilon]/\varepsilon^2)$ be the base space
of first
order infinitesimal deformations.   We turn our attention to a subfunctor
$\kd e\!f^\prime_{C,x} \subset \kd e\!f_{C,x}$ such that
\[
(T^1)^\prime := \kd e\!f_{C,x}^{\prime}(T_\varepsilon)
\]
is an ideal in $\kd e\!f_{C,x} (T_\varepsilon) \cong \C\{u,v\}/j(C,x)$ and the
corresponding
``global'' subfunctor $\kd e\!f^\prime_{C/S} \subset \kd e\!f_{C/S}$ where $\kd
e\!f^\prime_{C/S}(T)$
consists exactly of all those elements of $\kd e\!f_{C/S}(T)$ which are mapped
to
$\kd e\!f^\prime_{C,x}(T)$ for all points $x \in C$.
}
\end{sub}

{\bf Examples}:
    \begin{enumerate}
    \item
    $\kd e\!f^{ea}_{C/S}$ the subfunctor of $\kd e\!f_{C/S}$ consisting of all
    isomorphism classes of equianalytic deformations of $C/S$, that is of those
    deformations whose induced deformations of the analytic germs $(C,x)$
    happen to be equianalytic for all $x \in C$.   Here $\kd
e\!f^{ea}_{C,x}(T_\varepsilon) =
    0$ in  $\C\{u,v\}/j(C,x)$.

    \item
    $\kd e\!f^{es}_{C/S}$ the subfunctor of $\kd e\!f_{C/S}$ consisting of all
    isomorphism classes of equisingular deformations of $C/S$, that is of
    those deformations whose induced deformations of $(C,x)$ are
    equisingular for all $x \in C$.   Here $\kd e\!f^{es}_{C,x}(T_\varepsilon)
= I^{es}
    (C,x)/j(C,x)$.

    \item
Further examples are the equimultiple, equigeneric and equiclassical
deformation functors (cf.\ \cite{DH}).
    \end{enumerate}

\begin{remark}{\rm  $\kd e\!f^{ea}_{C/S}$ coincides with $\kd e\!f^{es}_{C/S}$
if
(and only if) $C$ has only simple (ADE)--singularities.
}
\end{remark}

Now, let $J_C$ be the ideal sheaf of $C$ in $\ko_S$, then we have the natural
exact sequence
\[
0 \to J_C/J^2_C \to \Omega^1_S \otimes_{\ko_S} \ko_C \to \Omega^1_C \to 0
\]
respectively its dual
\[
0 \to \theta_C \to \theta_S \otimes_{\ko_S} \ko_C \buildrel\Psi\over\to
\kn_{C/S} \to \kt^1_C \to 0.
\]

Here $\kn_{C/S} = \ko_S(C) \otimes_{\ko_S} \ko_C$ denotes the normal sheaf of
$C$ in $S$ and $\kt^1_C := \mbox{ Coker } (\Psi)$ is a skyscraper sheaf
concentrated in the singular points of $C$ with $H^0(C,\kt^1_C) \cong
\kd e\!f_C(T_\varepsilon)$ and with stalk in $x \in C$ equal to $\kt^1_{C,x}
\cong
\kd e\!f_{C,x}(T_\varepsilon) = T^1_{(C,x)}$ (for
 details cf.\ \cite{Art}).   Furthermore, for each subfunctor $\kd
e\!f^\prime_{C/S}$
as above, let $(\kt^1_C)'$ denote the subsheaf of $\kt^1_C$ with stalk in $x$
isomorphic to $(T^1)^\prime \subset T^1$ and
\[
\kn'_{C/S} := \mbox{ Ker}(\kn_{C/S} \to \kt^1_C/(\kt^1_C)').
\]
In particular,
\vspace{-1cm}

\begin{eqnarray*}
\kn^{ea}_{C/S} & = & \mbox{ Ker }(\kn_{C/S} \to \kt^1_C),\\
\kn^{es}_{C/S} & = & \mbox{ Ker } (\kn_{C/S} \to \kt^1_C/(\kt^1_C)^{es})
\end{eqnarray*}
where $\kt^1_{C,x} \cong \C \{u,v\}/j(C,x),\; (\kt^1_C)^{es}_x \cong
I^{es}(C,x)/j(C,x)$.

\begin{lemma}\label{3.2}  There is a canonical isomorphism
\[
\Phi: \quad \kd e\!f^\prime_{C/S}(T_\varepsilon)
\buildrel\cong\over\longrightarrow
H^0(C,\kn'_{C/S}).
\]
\end{lemma}

{\bf Proof}:  Each representative of an element in $\kd
e\!f^\prime_{C/S}(T_\varepsilon)$
is given by local equations $(f_i + \varepsilon g_i = 0)_{i\in I}$ (where
$f_i,g_i \in \Gamma(U_i,\ko_S)$ for an open covering $(U_i) _{i\in I}$ of
$S$), which satisfy
    \begin{itemize}
    \item $(f_i = 0)_{i\in I}$ are local equations for $C \subset S$

    \item $f_i + \varepsilon g_i = (a_{ij} + \varepsilon b_{ij}) \cdot (f_j +
    \varepsilon g_j)$ on $U_i \cap U_j =: U_{ij}$ with $a_{ij}$ a unit in
    $\Gamma(U_{ij}, \ko_S)$ and $b_{ij} \in \Gamma(U_{ij}, \ko_S)$

    \item the germ $g_{i,x}$ of $g_i$ projects to an element of
    $\kd e\!f^\prime_{C,x}(T_\varepsilon) = (T^1)^\prime$ for all $x \in C \cap
U_i$.
    \end{itemize}

For the induced sections $\frac{g_i}{f_i} \in \Gamma(U_i, \ko_S(C))$ it
follows immediately
\[
\frac{g_i}{f_i} - \frac{g_j}{f_j} = \frac{a_{ij}g_j + b_{ij}f_j}{a_{ij}f_j} -
\frac{g_j}{f_j} = \frac{b_{ij}}{a_{ij}} \equiv 0 \in \Gamma(U_{ij}, \kn_{C/S})
\]
and $\frac{g_i}{f_i}$ maps to an element of $(\kt^1_C)' \subset \kt^1_C$.
Hence,
$(\frac{g_i}{f_i})_{i \in I}$ defines a global section in $\kn'_{C/S}$.   It is
easy to check that in this way we get the isomorphism we were looking
for (cf.\ \cite{Mu}, \cite{Lo}).\hfill $\Box$

\begin{sub}\label{3.3}
{\rm
Let $C$ be a compact reduced curve, $\kf$ and $\kg$ coherent torsion--free
sheaves on $C$,
which have rank 1 on each irreducible component $C_i$ of $C$ and $x \in C$.
Then we define the {\sl local isomorphism defect} of $\kf$ in $\kg$ in $x$ as
\[
\mbox{isod}_x (\kf, \kg) := \min(\dim_{\C} \mbox{ Coker} (\varphi :
\kf_x \to \kg_x))
\]
where the minimum is taken over
all (injective) local homomorphisms $\varphi : \kf_x \to \kg_x$.   In
particular,
isod$_x(\kf,\kg)$ is a non--negative integer and not zero only in finitely many
points (in \cite{GrK} isod$_x(\kf,\kg)$ was denoted by $-$ind$_x(\kf,\kg)$).
We call
\[
\mbox{isod}(\kf,\kg) := \sum_{x \in C} \mbox{ isod}_x (\kf,\kg)
\]
the {\sl total (local) isomorphism defect} of $\kf$ in $\kg$.   For an
irreducible component
$C_i$ of $C$ and $x \in C_i$ we set
\[
\begin{array}{l}
\kf_{C_i} := \kf \otimes \ko_{C_i}\mbox{ modulo torsion}\\
\mbox{isod}_{C_i,x}(\kf,\kg) := \min(\dim_\C \mbox{ Coker} (\varphi_{C_i} :
\kf_{C_i,x} \to \kg_{C_i,x}))
\end{array}
\]
where the  minimum is taken over all $\varphi_{C_i}$, which are induced by
local homomorphisms $\varphi : \kf_x \to \kg_x$, and
\[
\mbox{isod}_{C_i}(\kf, \kg) := \sum_{x\in C_i} \mbox{ isod}_{C_i,x} (\kf, \kg).
\]

Note that this is again a non--negative integer.   In Chapter 4 we present
some explicit calculations.}
\end{sub}

\begin{proposition}\label{3.4}
(\cite{GrK}, Proposition 5.2)

Let $S$ be a smooth surface, $C \subset S$ a compact reduced curve and $\kf$ a
torsion--free coherent $\ko_C$--module  which
has rank 1 on each irreducible component $C_i$ of $C$  $(i = 1, \ldots, s)$.
Then $H^1(C,\kf) = 0$ if for $i = 1, \ldots, s$
\[
\chi(\kf_{C_i}) > \chi (w_{C,C_i}) - \mbox{ isod}_{C_i} (\kf,w_C).
\]

Here $\chi(\km) = \dim H^0((C,\km) - \dim H^1(C,\km)$ for a
coherent sheaf $\km$ on $C$ and  $w_C$ denotes the dualizing sheaf, $w_{C,C_i}
:= w_C \otimes
\ko_{C_i}$.
\end{proposition}

\begin{remark}{\rm  Using Riemann--Roch and the adjunction formula, the
condition above reads
\[
\deg(\kf_{C_i}) > (K_S + C) \cdot C_i - \mbox{ isod}_{C_i}(\kf,\ko_C)
\]
where $K_S$ is the canonical divisor on $S$.   Since isod is a local invariant
and since $C$ has planar singularities, we can replace $w_C$ by $\ko_C$.
}
\end{remark}

\begin{theorem}\label{3.5}
Let $S$ be a smooth complex surface and $C \subset S$ a reduced compact curve,
$H^{ea}_S$ respectively $H^{es}_S$ be the representing spaces for the
equianalytic respectively equisingular Hilbert functor, then
    \begin{itemize}
    \item[(i)] $\dim(H^{ea}_S, C) \ge C^2 + 1 - p_a(C) - \tau(C)$,
    with $\tau(C)
    = \sum_{x \in Sing(C)} \tau(C,x),$\\[1.0ex]
    $\dim(H^{es}_S, C) \ge C^2 + 1 - p_a(C) - \tau^{es}(C)$, with $\tau^{es}(C)
=
    \sum_{x \in Sing(C)} \tau^{es}(C,x),$\\[1.0ex]
    where $\tau(C,x) = \dim_\C(\ko_{C,x}/j(C,x))$ and $\tau^{es}(C,x) =
    \dim_\C(\ko_{C,x}/I^{es}(C,x))$

    \item[(ii)] If $H^1(C, \kn^{ea}_{C/S}) = 0$ (respectively
    $H^1(C,\kn^{es}_{C/S}) = 0)$ then $H^{ea}_S$ (respectively $H^{es}_S$) is
    smooth at $C$ of dimension
    \[
    C^2 + 1 - p_a (C) - \tau(C) \quad (\mbox{respectively } C^2 + 1 - p_a(C) -
    \tau^{es}(C)).
    \]

    \item[(iii)] Let $C = C_1 \cup \ldots \cup C_s$ be the decomposition into
    irreducible components, then
        \begin{itemize}
        \item[$\bullet$] $H^1(C, \kn^{ea}_{C/S}) = 0$ if for $i = 1,
        \ldots, s$
        \[
        - K_S \cdot C_i > D \cdot C_i + \tau(C_i) - \mbox{
        isod}_{C_i}(\kn^{ea}_{C/S}, \ko_C)
        \]

        \item[$\bullet$] $H^1(C,\kn^{es}_{C/S}) = 0$ if for $i = 1,
        \ldots, s$
        \[
        - K_S \cdot C_i > \sum\limits_{x \in\; Sing (C)} \dim_\C
        ((\ko_{C,x}/I^{es}(C,x)) \otimes \ko_{C_i,x}) - \mbox{ isod}_{C_i}
        (\kn^{es}_{C/S}, \ko_C)
        \]
        \end{itemize}

    where $D = \cup_{j\not= i} C_j$ and $K_S$ denotes the
    canonical divisor on $S$.   Moreover, the isomorphism defects
    isod$_{C_i}(\kn^{ea}_{C/S}, \ko_C)$ (respectively
isod$_{C_i}(\kn^{es}_{C/S},
    \ko_C))$ have the lower  bound $\#(C_i \cap \mbox{ Sing } C)$.
    \end{itemize}
\end{theorem}

\begin{remark}{\rm
\begin{enumerate}
\item 
If all singularities of $C$ are quasi--homogeneous or {\sl ordinary $k$--tuple
points}   (all branches are smooth with distinct tangents) then we obtain as an
equivalent criterium for the vanishing of $H^1(C,\kn^{es}_{C/S})$
\[
- - - K_S \cdot C_i > D \cdot C_i + \tau^{es}(C_i) - \mbox{ isod}_{C_i}
(\kn^{es}_{C/S},\ko_C).
\]

\item 
{}From the adjunction formula we obtain
\[
- - - K_S \cdot C_i = C^2_i - 2 p_a (C_i) + 2.
\]
\item 
If $C$ is irreducible, we have
    \begin{itemize}
    \item $H^1(C,\kn^{ea}_{C/S}) = 0$ if $-K_SC >  \tau(C) - \mbox{
    isod}(\kn^{ea}_{C/S}, \ko_C)$
    \item $H^1(C,\kn^{es}_{C/S}) = 0$ if $-K_SC >  \tau^{es}(C) - \mbox{
    isod}(\kn^{es}_{C/S}, \ko_C)$.
    \end{itemize}
\end{enumerate}}
\end{remark}

\begin{sub}\label{3.6}
{\rm {\bf Proof}:  Most parts of the proof are identical for the equianalytic
$(ea)$ and the equisingular $(es)$ case, there we use again the notation $H'$
respectively $\kn'_{C/S}$ as above:
\vspace{-0.5cm}

    \begin{itemize}
    \item[(ii)]
    Let $H^1(C, \kn'_{C/S}) = 0$ and $A \twoheadrightarrow A/(\eta) = \bar{A}$
    be a small extension of Artinian $\C$--algebras.   For the smoothness of
    $(H', C)$, we have to show that each equianalytic (respectively
    equisingular) family $\bar{\kc}$ over $\bar{A}$ lifts to an equianalytic
    (respectively equisingular) family $\kc$ over $A$.

    $\bar{\kc}$ is given by local equations $\bar{F_i} \in \Gamma(U_i, \ko_S
    \otimes \bar{A})$, where $(U_i)_{i \in I}$ is an open covering of $S$,
    such that

    $\bullet$ on $U_{ij} := U_i \cap U_j$, $\bar{F}_i = \bar{G}_{ij}
    \cdot \bar{F}_j$ with a unit $\bar{G}_{ij}$

    $\bullet$ the image $F^{(0)}_i$ of $\bar{F}_i$ in $\Gamma(U_i, \ko_S
    \otimes \C)$ is a local equation for $C \subset S$

    $\bullet$ the germs $\bar{F}_{i,x}$ describe an equianalytic
    (respectively equisingular) deformation of $(C,x)$.

    On the other hand, we know that the equianalytic (respectively
    equisingular) functor $E'$, which associates to each Artinian local
    $\C$--algebra the set of all equianalytic (respectively equisingular)
    deformations of $(C,x)$ over Spec $A$, is smooth and has a very good
    deformation theory (cf.\ \cite{Wa1},
    for the equisingular case).   Using the results of M.\
    Schlessinger (\cite{Schl}, Remark 2.17),
    this guarantees in particular

    \begin{itemize}
    \item[$\bullet$] the existence of an equianalytic (respectively
    equisingular) lifting $F_i \in \Gamma(U_i, \ko_S \otimes A)$ of
    $\bar{F}_i$

    \item[$\bullet$] for any lifting $G_{ij} \in \Gamma(U_{ij}, \ko_S \otimes
A)$
    of $\bar{G}_{ij}$ the existence of $h_{ij} \in \Gamma(U_{ij}, \ko_S \otimes
A)$
    with $F_i = G_{ij} \cdot F_j + \eta \cdot h_{ij}$ and $(h_{ij})_x \in
    j(C,x)$ (respectively $(h_{ij})_x \in I^{es}(C,x)$).
    \end{itemize}

    To obtain the lifted family we are looking for, we have to modify the
    $F_i$ and $G_{ij}$ in a suitable way, such that the $h_{ij}$ become 0.
    We know
    \begin{eqnarray*}
    \eta \cdot h_{ij} + \eta \cdot G_{ij} \cdot h_{jk} & = & F_i - G_{ij}
    \cdot F_j + G_{ij} \cdot (F_j - G_{jk} \cdot F_k)\\
    & = & \eta \cdot h_{ik} + (G_{ik} - G_{ij} \cdot G_{jk}) \cdot F_k
    \end{eqnarray*}
    where $(G_{ik} - G_{ij} \cdot G_{jk}) \in \Gamma(U_{ijk}, \ko_S \otimes
    (\eta))$ and $(\eta) \cdot \frak m_A = 0$.   As sections in $\ko_S
    \otimes A/\frak m_A = \ko_S \otimes \C$ we obtain
    \[
    h_{ij} + G_{ij}^{(0)} \cdot h_{jk} = h_{ik} + \left[\frac{1-G_{ij} \cdot
    G_{jk} \cdot G_{ik}^{-1}}{\eta}\right] \cdot G^{(0)}_{ik} \cdot F_k^{(0)}.
    \]

    Furthermore, $F_i^{(0)} = G_{ij}^{(0)} \cdot F_j^{(0)}$, which implies in
    $\Gamma(U_{ijk}, \kn_{C/S})$ the cocycle condition
    \[
    \frac{h_{ij}}{F_i^{(0)}} + \frac{h_{jk}}{F_j^{(0)}} =
    \frac{h_{ik}}{F_i^{(0)}}.
    \]

    From the definition of the $h_{ij}$ it follows that
    $\left(\frac{h_{ij}}{F_i^{(0)}} \mid i, j \in I\right)$ represents an
    element in $H^1(C,\kn'_{C/S}) = 0$.   Hence, there exist $f_i \in
    \Gamma(U_i, \ko_S \otimes \C)$ such that
    \[
    \frac{h_{ij}}{F_i^{(0)}} = \frac{f_j}{F_j^{(0)}} - \frac{f_i}{F_i^{(0)}}
    \]
    as sections in $\kn'_{C/S}$, especially $h_{ij} + f_i - f_j \cdot
    G_{ij}^{(0)} \in \Gamma(U_{ij}, J_C)$ and all germs $(f_i)_x$ lie in the
    Jacobian (respectively equisingularity) ideal.   Defining $g_{ij} :=
    \frac{h_{ij} + f_i - f_j \cdot G_{ij}^{(0)}}{F_j^{(0)}},\; \tilde{F}_i :=
    F_i + \eta \cdot f_i$ and $\widetilde{G}_{ij} := G_{ij} + \eta \cdot
    g_{ij}$ we obtain the lifted family we were looking for.

    \item[(i)] The germ $(H'_S, C)$ is the fibre over the origin of a
    (non--linear) obstruction map $H^0(C, \kn'_{C/S}) \to H^1(C, \kn'_{C/S})$
    (cf.\ \cite{La}, Theorem 4.2.4).   Hence
    \vspace{-0.5cm}

    \begin{eqnarray*}
    \dim\; H^0(C, \kn'_{C/S}) & \ge & \dim(H'_S, C)\\
    & \ge & \dim(H^0(C, \kn'_{C/S})) - \dim(H^1(C, \kn'_{C/S}))\\
    & = & \chi(\kn_{C/S}) - \chi(\kt^1_C/(\kt^1_C)')\\
    & = & \deg(\kn_{C/S}) + \chi(\ko_C) - \tau'(C)\\
    & = & C^2 + 1 -p_a(C) - \tau'(C)
    \end{eqnarray*}

    where $\tau'(C)$ denotes the total Tjurina number of $C$ (respectively
    $\tau^{es}(C))$.   Both inequalities become an  equality if
$H^1(C,\kn'_{C/S}) =
    0$.

    \item[(iii)] By Proposition \ref{3.4}, $H^1(C, \kn'_{C/S}) = 0$, if for $i
= 1,
    \ldots, s$
    \[
    \deg(\overline{\kn'_{C/S} \otimes \ko_{C_i}}) > (K_S + C) \cdot C_i -
    \mbox{ isod}_{C_i} (\kn'_{C/S}, \ko_C)
    \]
    where $\overline{\phantom{xxx}}$ denotes reduction modulo torsion.   On
    the other hand, we have the exact sequence
    \[
    0 \to \overline{\kn'_{C/S} \otimes \ko_{C_i}} \to \kn_{C/S} \otimes
    \ko_{C_i} \to (\kt^1_C/(\kt^1_C)^\prime), \otimes \ko_{C_i} \to 0
    \]
    which implies
    \[
    \deg(\overline{\kn'_{C/S} \otimes \ko_{C_i}}) = C \cdot C_i - \dim_\C
    H^0(C, (\kt^1_C/(\kt^1_C)^\prime) \otimes \ko_{C_i}).
    \]

    Finally, we obtain the above criteria by
    \[
    \dim_\C H^0(C, \kt^1_C/(\kt^1_C)^{es} \otimes \ko_{C_i}) = \sum_{x \in C}
    \dim_\C((\ko_{C,x}/I^{es}(C,x)) \otimes \ko_{C_i,x})
    \]
    respectively using the Leibniz rule (with $g$ as the equation of
    $(D,x))$ by
    \vspace{-0.5cm}

    \begin{eqnarray*}
    \dim_\C H^0(C, \kt^1_C \otimes \ko_{C_i}) & = & \sum_{x \in C} \dim_\C
    ((\ko_{C,x}/j(C,x)) \otimes \ko_{C_i,x})\\
    & = & \sum_{x \in C} \dim_\C(\ko_{C_i,x}/g \cdot j(C_i,x))\\
    & = & \sum_{x \in C} (\dim_\C(\ko_{C_i,x}/j(C_i,x)) +
    \dim_\C(\ko_{C_i,x}/g \cdot \ko_{C_i,x}))\\
    & = & \tau(C_i) + C_i \cdot D.
    \end{eqnarray*}
    \end{itemize}
\begin{flushright} $\Box$\end{flushright}
}
\end{sub}

\begin{sub}\label{3.7}{\rm
{\bf Curves in} $\P^2(\C)$

Let $C \subset \P^2 := \P^2(\C)$ be a reduced curve of degree $d$ with
(homogeneous) equation $F(X,Y,Z) = 0$, then we define the {\sl polar} of $C$
relative to
the point ($\alpha : \beta : \gamma) \in \P^2$ to be the curve
$C_{\alpha\beta\gamma}$ with equation $\alpha \cdot F_X(X,Y,Z) + \beta \cdot
F_Y(X,Y,Z) + \gamma \cdot F_Z(X,Y,Z) = 0$.
}\end{sub}

\begin{lemma}
The {\sl generic polar} $C_{\alpha\beta\gamma}$ (with
$(\alpha:\beta:\gamma) \in \P^2$ a generic point) is an irreducible curve of
degree $d-1$, if and only if $C$ is not the union of $d \ge 3$ lines through
the same point.
\end{lemma}

{\bf Proof}:  Applying Bertini's theorem, we have irreducibility if there is no
algebraic relation between $F_X, F_Y$ and $F_Z$.   Considering such a
(homogeneous) relation of minimal degree
\[
\sum_\alpha a_\alpha F^{\alpha_1}_X \cdot F^{\alpha_2}_Y \cdot F^{\alpha_3}_Z =
0
\]
and differentiating, we obtain a system of equations
\[
A \cdot \;\left(\begin{array}{c}F_{XX}\\F_{XY}\\F_{XZ}\end{array}\right)
\; + B \cdot\;\left(\begin{array}{c}F_{YX}\\F_{YY}\\F_{YZ}\end{array}\right)
\; + \Gamma
\cdot\;\left(\begin{array}{c}F_{ZX}\\F_{ZY}\\F_{ZZ}\end{array}\right)
\; = 0
\]
where $A, B$ and $\Gamma$ generically do not vanish.   Now the lemma follows
from the fact that the Hessian covariant vanishes identically only if $C$ is
the union of $d$ lines through one point (cf.\ \cite{He}, Lehrsatz 6).
\hfill $\Box$\\

In the following we choose suitable coordinates such that Sing $C$ lies in the
affine plane $Z \not= 0$ and a
generic polar $C' \subset P^2$ relative to $(\alpha : \beta : 0)$ with equation
$\alpha \cdot F_X + \beta
\cdot F_Y = 0$ is irreducible.

Then we have an obvious morphism $\ko_{C'}(d) \to \kt^1_C$ given by the natural
projections
\[
\ko_{C',x} = \ko_{\P^2,x}/(\alpha f_X + \beta f_Y) \to \ko_{\P^2,x}/j(C,x) =
\kt^1_{C,x}
\]
where $f(X,Y) = F(X,Y, 1)$  is the affine equation of $C$.

We define
\vspace{-0.5cm}

\begin{eqnarray*}
\tilde{\kn}^{ea}_{C'/\P^2} & := & Ker(\ko_{C'}(d) \to \kt^1_C)\\
\tilde{\kn}^{es}_{C'/\P^2} & := & Ker(\ko_{C'}(d) \to \kt^1_C/(\kt_C^1)^{es}).
\end{eqnarray*}

\begin{corollary}\label{3.8}
Let $C \subset \P^2$ be a reduced projective curve of degree $d$,
$C_i(i = 1,\ldots, s)$ its irreducible components and $d_i$
the degree of $C_i$.

    \begin{enumerate}
    \item[(i)] $H^1(C, \kn^{ea}_{C/\P^2}) = 0$ if and only if the forgetful
    morphism $\kd e\!f_{C/\P^2} \to \prod_{x \in Sing\, C}$ $\kd e\!f_{C,x}$ is
    surjective.

    \item[(ii)] If $H^1(C, \kn^{ea}_{C/\P^2}) = 0$ (respectively $H^1(C,
    \kn^{es}_{C/S}) = 0)$ then $(H^{ea}_{\P^2}, C)$ (respectively
    $(H^{es}_{\P^2},C)$)
    is smooth of dimension $\frac{1}{2} d(d+3) - \tau(C)$ (respectively
    $\frac{1}{2} d \cdot (d+3) - \tau^{es}(C))$.

    \item[(iii)] $H^1(C, \kn^{ea}_{C/\P^2}) = 0$ if for $i = 1, \ldots, s$
    \[
    3d_i > d_i \cdot (d - d_i) + \tau(C_i) - \mbox{ isod}_{C_i}
    (\kn^{ea}_{C/\P^2}, \ko_C),
    \]
moreover, isod$_{C_i}(\kn^{ea}_{C/\P^2}, \ko_C) \ge \#$ (Sing$(C) \cap C_i)$,\\
    $H^1 (C,\kn^{es}_{C/S}) = 0$ if for $i = 1, \ldots, s$
    \[
    3 d_i > d_i \cdot (d-d_i) + \tau^{es}(C_i) - \mbox{
    isod}_{C_i}(\kn^{es}_{C/\P^2}, \ko_C).
    \]

    \item[(iv)] If $C$ is not the union of $d \ge 3$ lines through one point
    $H^1(C, \kn^{ea}_{C/\P^2})$ (respectively $H^1(C, \kn^{es}_{C/\P^2})$)
vanishes,
    if
    \[
    4 d > 4 + \tau(C) - \mbox{ isod}(\tilde{\kn}^{ea}_{C'/\P^2}, \ko_{C'})
    \]
    \[
    (\mbox{respectively } 4 d > 4 + \tau^{es}(C) - \mbox{ isod}
    (\tilde{\kn}^{es}_{C'/\P^2}, \ko_{C'}))
    \]
    where $C'$ denotes the generic polar as defined above.
    \end{enumerate}
\end{corollary}

{\bf Proof}:  (ii) and (iii) follow immediately from Theorem \ref{3.5}.   To
prove (iv), we consider the exact sequences
\[
\begin{array}{ccccccccc}
0 & \to & \kn'_{C/\P^2} & \to & \kn_{C/\P^2} & \to & \kt^1_C/(\kt^1_C)' & \to &
0\\
0 & \to & \tilde{\kn}'_{C'/\P^2} & \to & \ko_{C'}(d) & \to & \kt^1_C/(\kt^1_C)'
&
\to & 0
\end{array}
\]
and the corresponding long exact cohomology sequences where $'$ represents
again both the equianalytic and the equisingular case.   We know that
$\kn_{C/\P^2} \cong \ko_C(d)$ and (by Proposition \ref{3.4}) $H^1(C,\ko_C(d))
= 0$.   Furthermore, if $C$ is not the union of $d \ge 3$ lines through one
point, $C'$ is
irreducible and
\vspace{-0.5cm}

\begin{eqnarray*}
\deg(\tilde{\kn}'_{C'/\P^2}) - (K_{\P^2} + C') \cdot C' & = &
\deg(\ko_{C'}(d)) - \tau'(C) - (d-4) \cdot (d-1)\\
& = & 4 \cdot (d-1) - \tau'(C).
\end{eqnarray*}

Hence, applying Proposition \ref{3.4} the conditions in (iv) guarantee the
vanishing of $H^1(C', \tilde{\kn}'_{C'/\P^2})$.   Additionally, the exact
sequence
\[
0 \to \ko_{\P_2} \to \ko_{\P_2} (d) \to \ko_C (d) \to 0
\]
respectively an analogous sequence for $C'$ induce surjective mappings
\[
\Phi : H^0(\P^2, \ko_{\P^2}(d)) \twoheadrightarrow H^0(C, \ko_C(d)) \mbox{
respectively }
\Phi' : H^0(\P^2, \ko_{\P^2}(d))\twoheadrightarrow H^0(C', \ko_{C'}(d))
\]
which lead to a commutative diagram with exact horizontal rows

\unitlength1cm
\begin{picture}(10,3)
\put(0.5,1.5){$H^0(\P^2, \ko_{\P^2}(d))$}
\put(3.5,1.8){\vector(2,1){1.2}}
\put(3.8,2.2){$\Phi$}
\put(5.0,2.4){$H^0(C, \ko_C(d)) \to H^0(C, \kt^1_C/(\kt^1_C)')
\to H^1(C, \kn'_{C/\P^2}) \to 0$}
\put(9.0,1.5){$\|$}
\put(3.8,0.75){$\Phi'$}
\put(3.5,1.4){\vector(2,-1){1.2}}
\put(5.0,0.5){$H^0(C', \ko_{C'}(d)) \to H^0(C, \kt^1_C/(\kt^1_C)') \to 0.$}
\end{picture}

This shows that $H^0(C,\ko_C(d)) \to H^0(C,\kt^1_C/(\kt^1_C)')$ is surjective
and, hence, $H^1(C, \kn'_{C/\P^2}) = 0$.

The same argument shows that $H^0(C, \ko_{C}(d)) \to H^0(C,\kt^1_C)$ is
surjective if and only if $H^1(C, \kn^{ea}_{C/\P^2}) = 0$.   Since
$\kd e\!f_{C/\P^2}(T_\varepsilon) = H^0(C, \ko_C(d))$
and $H^0(C, \kt^1_C) \cong \prod_{x\in Sing\, C}$ $\kd
e\!f_{C,x}(T_\varepsilon)$
this is equivalent to $\kd e\!f_{C/\P^2}(T_\varepsilon) \to \prod_{x
\in Sing(C)}$
$\kd e\!f_{C,x}(T_\varepsilon)$ being surjective (and $\kd e\!f_{C/\P^2}$ being
unobstructed).   But $\kd e\!f_{C/\P^2}$ and $\prod$ $\kd e\!f_{C,x}$ being
unobstructed, the surjectivity on the tangent level implies the surjectivity
of the functors.\hfill $\Box$

\begin{remark}{\rm
\begin{itemize}
\item[(i)] We call the inequalities in \ref{3.8} (iii) respectively \ref{3.8}
(iv) the $3d$-- respectively $4d$--criteria.

\item[(ii)] The use of a generic polar is due to Shustin \cite{Sh1}, who
obtained (with a different proof) the weaker inequality $4d > 4 + \mu(C)$
instead of \ref{3.8} (iv), where $\mu(C)$ is the total Milnor number of $C$.
\end{itemize}
}
\end{remark}

\begin{sub}\label{3.9}{\rm {\bf A generalization}:
We are also interested in families of curves in $\P^2$ of degree $d$ where for
some singularities the analytic type is fixed, for others only the topological
type is fixed and for the remaining singularities any deformation is allowed.

Let $C \subset \P^2$ be of degree $d$ and Sing$(C) = \{x_1, \ldots, x_k\} \cup
\{y_1, \ldots, y_\ell\} \cup \{z_1, \ldots, z_m\}$.   We define the subsheaf
$(\kt^1_C)^\prime$ of $\kt^1_C$ by
\[
(\kt^1_C)'_x = \left\{
\begin{array}{ll}
0 & \mbox{ if } x \in \{x_1, \ldots, x_k\}\\
I^{es}(C,x)/j(C,x) & \mbox{ if } x \in \{y_1, \ldots, y_\ell\}\\
T^1_{(C,x)} & \mbox{else}
\end{array}\right.
\]
and put
\[
\tau'(C) := \dim_\C H^0(C,\kt^1_C/(\kt^1_C)^\prime) = \sum\limits^k_
{i=1} \tau(C,x_i)
+ \sum\limits^\ell_{j=1} \tau^{es}(C,y_j).
\]
Assume there exists a reduced curve $C' \subset \P^2$ of degree $d'$ with the
following
properties:

\begin{itemize}
\item[(a)] $C'$ is irreducible,
\item[(b)] $\{x_1, \ldots, x_k\} \cup \{y_1, \ldots, y_\ell\} \subset C'$,
\item[(c)] if $f_j$ is a local equation of $(C', x_j),\; j = 1, \ldots, k$,
then $f_j \in j(C,x_j)$;  if $f_j$ is a local equation of $(C',y_j),\; j = 1,
\ldots, \ell$, then $f_j \in I^{es}(C,y_j)$.
\end{itemize}

Define
\vspace{-0.5cm}

\begin{eqnarray*}
\kn' & := & \mbox{Ker} (\kn_{C/\P^2} = \ko_C(d) \to \kt^1_C/(\kt^1_C)
^\prime),\\
\widetilde{\kn}' & := & \mbox{Ker} (\ko_{C'}(d) \to \kt^1_C/(\kt^1_C)^\prime).
\end{eqnarray*}

Let $\ka$ be the analytic singularity type defined by $(C,x_1), \ldots,
(C,x_k), \kt$ the topological singularity type defined by $(C,y_1), \ldots,
(C,y_\ell)$ and let $\kh ilb^{\ka,\kt}_{\P^2}$ denote the functor
parametrising
proper and flat families of reduced curves in $\P^2$ which have $k$ singular
points of fixed analytic type $\ka$ and $\ell$ singular points of fixed
topological type $\kt$ (see \ref{2.2} for a precise definition).   This functor
is represented by a locally closed subspace $H^{\ka,\kt}_{\P^2} \subset
H_{\P^2}$ (cf.\ Proposition 2.3).
}
\end{sub}

\begin{proposition}
Let $C \subset \P^2$ be a reduced projective curve of degree $d$, Sing$(C) =
\{x_1, \ldots, x_k\} \cup \{y_1, \ldots, y_\ell\} \cup \{z_1, \ldots,
z_m\}$
and assume that there exists a curve $C' \subset \P^2$ of degree $d'$
satisfying (a) -- (c) above.
\vspace{-0.5cm}

\begin{itemize}
\item[(i)] If $\ell = 0$, then $H^1(C,\kn') = 0$ if and only if $\kd
e\!f_{C/\P^2}
\to \prod\limits^k_{i=1}\kd e\!f_{(C,x_i)}$ is surjective.

\item[(ii)] If $H^1(C,\kn') = 0$ then $H^{\ka,\kt}_{\P^2}$ is smooth at $C$ of
dimension $\frac{1}{2} d (d+3) - \tau'(C)$.

\item[(iii)] $H^1(C,\kn') = 0$ if $d'(d-d'+3) > \tau'(C) -$
isod$(\widetilde{\kn}',\ko_{C'})$.
\end{itemize}
\end{proposition}

{\bf Proof}:  Use the same argumentation as for Corollary 3.12.\hfill $\Box$

\begin{remark}{\rm
\begin{enumerate}
\item If $\ell = m = 0$, (i) is the same as 3.12 (i).   If $\ell = m = 0$
(respectively $k = m = 0)$ (ii) is the same as 3.12\ (ii).   If $C$ is
irreducible we may take $C' = C$ and then (iii) is equivalent to 3.12 (iii)
for $\ell = m = 0$ respectively $k = m = 0$.   If $C$ is not the union of $d$
lines through one point, we may take $C'$ to be a generic polar and then (iii)
is the same as 3.12 (iv) for $\ell = m = 0$ (respectively $k = m = 0$).

\item We obtain the best possible result for a curve $C'$ of degree $d' =
\frac{d+3}{2}$ satisfying (a) -- (c).   In \cite{Sh3}, Shustin has proven the
existence of such an irreducible curve $C'$ of degree
\[
d' \le (2\kappa^2 + \sqrt{\kappa}) \sqrt{\mu(C)} + (1 - \frac{1}{\kappa}) d,
\]
where $\mu(C)$ denotes the total Milnor number and
\[
\begin{array}{ll}
\kappa = &\; \max\;\{\mu(x_i, C) + \mbox{ mult}_{x_i}(C),\; \mu(y_j,C) + \mbox{
     mult}_{y_j}(C)\} - 1\\[-0.5ex]
    & {1 \le i \le k\atop  1 \le j \le \ell}
\end{array}
\]
\end{enumerate}
}
\end{remark}
\newpage

\section{Local isomorphism defects of plane curve singularities}

Let $u, v$ be local coordinates of the smooth surface $S$ in a singular
point $x$ of the reduced compact curve $C \subset S,\; C = C_1 \cup \ldots
\cup C_s$ the decomposition into irreducible components and $f(u,v) = 0$
(respectively $f_i(u,v) = 0$) be local equations of $(C,x)$ (respectively
($C_i,x$)).   In the following, we give estimations for the (local)
isomorphism defects occurring in Chapter 3:

\begin{lemma}\label{4.1}
For a  reduced plane curve singularity $(C,x) \subset (S,x)$ we have:

\begin{itemize}
\item [(i)] isod$_x (\kn^{ea}_{C/S}, \ko_C) = 1$   if $(C,x)$ is
quasihomogeneous.\\
           isod$_x (\kn^{ea}_{C/S}, \ko_C) > 1$ if $(C,x)$ is
           not quasihomogeneous.

\item[(ii)] isod$_{C_i,x} (\kn^{ea}_{C/S}, \ko_C) \ge 1$ for $i = 1, \ldots,
s$.
\end{itemize}
\end{lemma}

{\bf Proof}:
\vspace{-0.5cm}

\begin{itemize}
\item[(i)] By definition isod$_x (\kn^{ea}_{C/S}, \ko_C) = \min\;
\dim_\C(\mbox{
coker }\, \varphi : j(C,x) \to \ko_{C,x})$ where the minimum is taken over all
$\ko_{C,x}$--linear maps.   Now the Jacobian ideal $j(C,x)$ of an isolated
singularity cannot be generated by a single element and there is an
isomorphism $\varphi : j(C,x) \buildrel \cong\over\to \frak{m}_{C,x}$ exactly
if
$(C,x)$ is quasihomogeneous.

\item[(ii)] We have to look for an $\ko_{C,x}$--linear map $\varphi : j(C,x)
\to \ko_{C,x}$ whose restriction to $(C_i,x)$ has minimal cokernel.   The
above statement follows immediately. \hfill $\Box$
\end{itemize}

\begin{sub}\label{4.2}{\rm
Let $(C,x)$ be quasihomogeneous with positive weight vector $w = (w_1, w_2)$
and (weighted)  degree $d$, then (as an
$\ko_{C,x}$--ideal) $I^{es}(C,x)$ is generated by the Jacobian ideal $j(C,x)$
and all monomials $u^\alpha v^\beta$ with $w_1 \cdot \alpha + w_2 \cdot
\beta \ge d$.   Furthermore we have the normalization
\[
n : \ko_{C,x} \hookrightarrow \bar{\ko} := \prod\limits^r_{i=1} \C\{t\}
\]
where $r$ denotes the number of local irreducible components $(C^{(i)},x)$ of
$(C,x)$ (not to be confused with the global components $C_i$).

In the following, we use the notations cond$(\ko)$, cond$(j)$ respectively
cond$(I^{es})$ for the {\sl conductor ideals} of $\ko_{C,x}$, the Jacobian
respectively the equisingularity ideal in $\ko_{C,x}$, where, for an
$\ko_{C,x}$--ideal $I$,
\[
\mbox{ cond} (I) := \{ g \in I \mid g \cdot \bar{\ko} \subset I\}.
\]

Furthermore, for all these $\ko_{C,x}$--ideals, we denote by $\Gamma(I)
\subset \N^r$ the {\sl set of values} of $I$ and by $\underline{c} (I) \in
\N^r$ the
{\sl conductor} of $I$, that is $\Gamma(\mbox{cond}(I)) = \underline{c} (I) +
\N^r$.}
\end{sub}

\begin{lemma}
Let $(C,x)$ be quasihomogeneous of degree $d$,
then
\begin{itemize}
\item[(i)] isod$_x (\kn^{es}_{C/S}, \ko_C) = \delta(C,x) -
\dim_\C(I^{es}(C,x)/$cond$(I^{es}))$\\
{\it especially}:  isod$_x(\kn^{es}_{C/S}, \ko_C) \ge 1$ with equality if
and only if $j(C,x) = I^{es}(C,x)$.

\item[(ii)] isod$_{C_i,x} (\kn^{es}_{C/S}, \ko_C) =
\dim_\C((\ko_{C,x}/\mbox{cond}(\ko)) \otimes \ko_{C_i,x}) -
\dim_\C((I^{es}(C,x)/\mbox{cond}(I^{es})) \otimes \ko_{C_i,x})$.
\end{itemize}
\end{lemma}
\newpage

{\bf Proof}:
\begin{itemize}
\item[(i)] To calculate isod$_x(\kn^{es}_{C/S}, \ko_C)$ we have to consider an
$\ko_{C,x}$--linear mapping
\[
\Psi : I^{es}(C,x) \to \ko_{C,x}
\]
with minimal cokernel.   We know that such a $\Psi$ maps cond$(I^{es})$ to
cond$(\ko)$, hence we obtain the estimate
\vspace{-0.5cm}

\begin{eqnarray*}
\mbox{isod}_x(\kn^{es}_{C/S}, \ko_C) & \ge &
\dim_\C(\ko_{C,x}/\mbox{cond}(\ko))
- - - \dim_C(I^{es}(C,x)/\mbox{cond}(I^{es}))\\
& = & \delta(C,x) - \dim_\C(I^{es}(C,x)/\mbox{cond}(I^{es}))
\end{eqnarray*}
with equality if and only if there exists a $\Psi$ that maps cond$(I^{es})$
onto cond$(\ko)$, or equivalently (using $\varphi : j(C,x) \buildrel
\cong\over\to \frak m_{C,x})$ if and only if we can find an
$\ko_{C,x}$--linear mapping
\[
\Phi : I^{es}(C,x) \to j(C,x)
\]
of weighted degree $\underline{c} (j) - \underline{c} (I^{es})$.   Now
$\underline{c}(I^{es})-\underline{1}$ is a maximal (in the sense of
\cite{De}) in the semigroup $\Gamma(I^{es}) \supset \Gamma(j)$.
Hence, $(\underline{c}(\ko) - \underline{c}(j)) + \underline{c} (I^{es}) -
\underline{1}$
is a maximal in $\Gamma(\ko)$ and using the symmetry of $\Gamma(\ko)$ we see
that
\[
(\underline{c}(\ko) - \underline{1}) - (\underline{c}(\ko) - \underline{c}(j)
+ \underline{c}
(I^{es}) - \underline{1}) = \underline{c}(j) - \underline{c}(I^{es}) \in
\Gamma(\ko).
\]

The additional statement is an immediate consequence from the fact that all
monomials of degree at least $d$ are contained in
cond$(I^{es})$;  thus,
\vspace{-0.5cm}

\begin{eqnarray*}
\dim_\C(I^{es}(C,x)/\mbox{cond}(I^{es})) & = &
\dim_\C(j(C,x)/\mbox{cond}(I^{es}) \cap j(C,x))\\
& \le & \dim_\C(j(C,x)/\mbox{cond}(j))\\
& = & \delta(C,x) - 1.
\end{eqnarray*}

\item[(ii)] Follows from the considerations above.\hfill $\Box$
\end{itemize}

\begin{remark}{\rm
  If $(C,x)$ is quasihomogeneous, then it is
easy to see that
\vspace{-0.5cm}

\begin{eqnarray*}
\dim_\C((\ko_{C,x}/\mbox{cond}(\ko)) \otimes \ko_{C_i,x}) & = &
\dim_\C(\ko_{C_i,x}/\mbox{cond}(\ko_{C_i})) + (D \cdot C_i,x)\\
\dim_\C((I^{es}(C,x)/\mbox{cond}(I^{es})) \otimes \ko_{C_i,x}) & \ge &
\dim_\C(I^{es}(C_i,x)/\mbox{cond}(I^{es}(C_i,x)))
\end{eqnarray*}

(where $C_j$, $j = 1, \ldots, s)$ are the irreducible components of $C$ and
($D = \bigcup\limits_{j\not= i} C_j)$.   Thus we obtain as an upper bound
\[
\mbox{isod}_{C_i,x}(\kn^{es}_{C/S}, \ko_C) \le \mbox{ isod}_x(\kn^{es}_{C_i/S},
\ko_{C_i}) + (D \cdot C_i,x).
\]
Furthermore, for integer weights $w_1 \ge w_2$, gcd$(w_1, w_2) = 1$,  the
difference is bounded by $\frac{w_1 - 1}{w_2} + 2$.
}
\end{remark}

\begin{sub}\label{4.3}{\rm
{\bf Examples}:
\begin{itemize}
\item[(i)] If $(C,x)$ is an ADE--singularity, then isod$_x(\kn^{es}_{C/S},
\ko_C) = 1$.

\item[(ii)] If $(C,x)$ is homogeneous of degree $r \ge 3$, then
\[
\mbox{isod}_x(\kn^{es}_{C/S}, \ko_C) = \frac{r \cdot(r-1)}{2} -2 =
\tau^{es}(C,x) - \mbox{ mult}_x(C)
\]
furthermore, let $C = C_i \cup D$ as above, then we obtain for $(C_i,x)$ a
smooth branch
\[
\mbox{isod}_{C_i,x} (\kn^{es}_{C/S}, \ko_C) = r-2,
\]
while in the singular case
\[
\mbox{isod}_{C_i,x} (\kn^{es}_{C/S}, \ko_C) = (C_i \cdot D, x) + \mbox{ isod}_x
(\kn^{es}_{C_i/S}, \ko_{C_i}).
\]
More generally, these statements are valid, if $(C,x)$ is an ordinary
$r$--tuple point ($r$ smooth branches with different tangents, $r \ge 3$);
in this case each equimultiple deformation is equisingular.

\item[(iii)] If $(C,x)$ has the local equation $u^p - v^q = 0$ where $q \ge
p \ge 3$ and $(C_i,x)$ consist of $b \le r = \mbox{ gcd } (p,q)$ irreducible
branches, then
\[
\mbox{isod}_{C_i,x} (\kn^{es}_{C/S}, \ko_C) = \frac{b}{2r} \cdot
(pq(2-\frac{b}{r}) + (r-p-q)) - \Big[\frac{q-2}{p}\Big] - M,
\]
where $M = 2$ unless $b = 1$ and $q = k \cdot p\;\; (k \in \N)$, then $M = 1$.
\end{itemize}
}
\end{sub}

\begin{sub}\label{4.4}{\rm
Now let $S = \P^2$ and $C \subset \P^2$ be different from the union of $d \ge
3$ lines through one point and $C' \subset \P^2$ denote the (irreducible)
generic polar of $C$ (cf.\ \ref{3.7}) with affine equation $\alpha \cdot f_X
+ \beta f_Y = 0$.
}
\end{sub}

\begin{lemma}
\begin{itemize}
\item[(i)] isod$_x (\widetilde{\kn}^{ea}_{C'/\P^2}, \ko_{C'}) \ge 0$ with
equality if and only if $(C,x)$ is quasihomogeneous.
\item[(ii)] isod$_x (\widetilde{\kn}^{es}_{C'/\P^2}, \ko_{C'}) \ge \delta(C',x)
- - - \dim_\C(I^{es}(C,x) \otimes \ko_{C',x}/\mbox{cond}(I^{es}(C,x) \otimes
\ko_{C',x}))$.
\end{itemize}
\end{lemma}

{\bf Proof}:  By definition, equality in (i) is equivalent to the statement
that the $\ko_{C',x}$--ideal generated by $f, f_X$ and $f_Y$ is generated by
one single element.   Obviously, this holds exactly if $(C,x)$ is
quasihomogeneous, (ii) follows from the considerations in the proof of Lemma
4.3.\hfill
$\Box$

\begin{remark}{\rm  The generic polar $C'$ depends on the whole curve $C$ and
not only on the germ $(C,x)$, hence, in general it is {\sl not} enough to know
the local
equation of $(C,x)$ to be able to calculate
isod$_x(\widetilde{\kn}^{es}_{C'/\P^2}, \ko_{C'})$.   For example, if $(C,x)$
is
homogeneous of degree $d \ge 6$, $(C',x)$ need not be quasihomogeneous.   But,
in special cases, we are able to give explicit formulas:
}
\end{remark}

\begin{sub}\label{4.5}{\rm
{\bf Examples}:

\begin{itemize}
\item[(i)] If $(C,x)$ is an ADE--singularity, then
isod$_x(\widetilde{\kn}^{es}_{C'/\P^2}, \ko_{C'}) = 0$.
\item[(ii)] If $(C,x)$ has the (homogeneous) local equation $(u^r - v^r = 0)\;
(r \ge 3)$
then $(C',x)$ has an equation $(\tilde{u}^{r-1} - \tilde{v}^{r-1} = 0)$ and
\[
\mbox{isod}_x (\widetilde{\kn}^{es}_{C'/\P^2}, \ko_{C'}) = \frac{r \cdot
(r-3)}{2}.
\]
Moreover,  in
the case of a (not homogeneous) local equation $(u^p - v^q = 0)\; (q > p \ge
3)$, $(C',x)$ has an equation $(\tilde{u}^{p-1} - \tilde{v}^{q-1} = 0)$ and we
obtain the estimate
\[
\mbox{isod}_x(\widetilde{\kn}^{es}_{C'/\P^2}, \ko_{C'}) \ge
\frac{(p-3)(q-1)}{2} + \frac{2gcd(p-1, q-1) - gcd(p,q) -1}{2}-
\Big[\frac{q}{p}\Big] + \varepsilon
\]
where $\varepsilon = 0$ unless $p$ divides $q$, then $\varepsilon = 1$.
\end{itemize}
}
\end{sub}

\begin{sub}{\rm {\bf Problem}:  Do the different isomorphism defects
considered above behave (lower)
semicontinuous under equianalytic respectively equisingular deformations of
$C/S$?
}
\end{sub}
\newpage

\section{Applications and Examples}
\begin{corollary}\label{5.1}
\begin{itemize}
\item[(i)]
Let $C \subset \P^2$ be a curve of degree $d$ with {\sl ordinary} $(
k_i${--}) {\sl multiple points} $(i = 1, \ldots, N)$ as the only
singularities, $C = C_1 \cup \ldots \cup C_s$ its decomposition into
irreducible components (deg $C_i = d_i)$ then
\[
(H^{es}_{\P^2}, C) \mbox{ is smooth of dimension }\frac{d
\cdot(d+3)}{2} - \sum\limits^N_{j=1} \bigg(\frac{k_j \cdot(k_j+1)}{2} -
2\bigg), \mbox{
if for } i = 1, \ldots, s
\]
\[
3 \cdot d_i > \sum_{{x \in C_i \cap Sing\, C\atop mult_x(C)>2}}
\mbox{mult}_x(C_i)
\]
where mult$_x(C)$ (respectively mult$_x (C_i)$) denote the multiplicity of $C$
(respectively $C_i$) at $x$.

\item[(ii)] Let $C \subset \P^2$ be a curve of degree $d$ whose singularities
are all of local equations $(u^{p_i} - v^{q_i} = 0)\; (q_i \ge p_i)$ or
ADE--singularities, then
\[
 (H^{es}_{\P^2}, C)\mbox{ is smooth of dimension }
\frac{d(d+3)}{2} - \sum\limits^N_{i=1} \bigg(\frac{(p_i + 1) \cdot (q_i + 1) -
gcd(p_i, q_i) - 5}{2} - \Big[\frac{q_i}{p_i}\Big] + \varepsilon_i\bigg)
\]
 where
$\varepsilon_i = 0$ unless $p_i$ divides $q_i$, then $\varepsilon_i = 1$, if
\[
4d > 4 + \sum_{\{ADE\}} \mu(C,x) + \sum_{\{not\; ADE\}} (p_i +
2q_i -3 - gcd(p_i -1, q_i - 1)).\]

\end{itemize}
\end{corollary}

{\bf Proof}:  The statements follow immediately from Corollary 3.12 and
Example 4.5 (ii), respectively Example 4.9 (ii). \hfill $\Box$

\begin{remark}\label{5.2}{\rm  The result in (i) was already obtained by C.\
Giacinti-Diebolt
(\cite{Gia}) using vanishing theorems on the normalization of $C$.   It implies
the ancient result of Severi \cite{Sev} that for a curve $C$ with no other
singularities but ordinary double points $(H^{es}_{\P^2}, C)$ is smooth.

Another consequence of the calculations in Chapter 4 is the following:  the
contribution of a quasihomogeneous singularity $(C,x)$ with local equation
$(u^p - v^q + uv \cdot \tilde{f}(u,v) = 0)$, $q \ge p \ge 3$, to the
right--hand side in the $3d$--criterion for an irreducible curve $C$ is
\[
\tau^{es}(C,x) - \mbox{ isod}_x(\kn^{es}_{C/S}, \ko_C) = p + q - gcd(p,q) -
\varepsilon
\]
where $\varepsilon = 0$, unless $q \equiv 1 \mbox{ mod } p$, then $\varepsilon
=
1$, while the contribution of an $A_k$--singularity is $k-1$.   This
corresponds to the result of E.\ Shustin in \cite{Sh2}.

Nevertheless, in some cases the new $4d$--criterion gives more information:
}
\end{remark}

\begin{sub}{\rm
{\bf Example}:
\begin{itemize}
\item[(a)] $(x^4 - x^2 z^2 + y^2 z^2 + y^3 z) \cdot y \cdot (x + 2y + z) \cdot
(x-2y-z) = 0$ defines a {\it reducible}\/ curve $C \subset \P^2$ having exactly
3 ordinary triple points lying on 1 line (hence, Corollary \ref{5.1} does {\it
not}\/ apply) and 7 ordinary double points.   But $4d = 28 > 23 = 4 +
\tau(C)$; hence, $(H^{ea}_{\P^2}, C)$, respectively $(H^{es}_{\P^2}, C)$ are
smooth of dimension 16.
\end{itemize}}
\end{sub}

\begin{sub}{\rm
In general, it is a difficult problem to determine for a given $d$ whether
there exists a projective plane curve of degree $d$ having a fixed number of
singularities of given {\sl analytic} type.   On the other hand, the local
deformations of a plane curve singularity are well understood.   Hence,
knowing about the existence of one low degree curve with ``big'' singularities
our $4d$--criterion allows us to give positive answers to some of the above
existing problems (answers which we did not obtain with the $3d$--criterion
in \cite{GrK}).}
\end{sub}

\begin{remark}{\rm
\begin{enumerate}
\item The surjectivity statement in \ref{3.8} (i) implies:  let $C \subset
\P^2$ be of degree $d$ such that $H^1(C,\kn^{ea}_{C/\P^2}) = 0$ and let $\{x_1,
\ldots, x_n\}$ be any subset of Sing$(C)$.   If, for $i = 1, \ldots, n$, the
germ $(C,x_i)$ admits a deformation with nearby fibre having singularities
$y^1_i, \ldots, y^{s_i}_i$, then $C$ admits an embedded deformation with
nearby curve $C_t \subset \P^2$ having $y^1_1, \ldots, y^{s_1}_1, \ldots,
y^1_n, \ldots, y^{s_n}_n$ as singularities.   Hence, there exists a curve of
degree $d$ with the $y$'s as singularities.

\item The $4d$--criterion in \ref{3.8} has the advantage that $C$ need not be
irreducible.   On the other hand, in the $3d$--criterion in \ref{3.5} and
\ref{3.8} we
can completely forget about $A_1$--singularities on $C$.   By Lemma 4.1,
for a node we have $\tau(C_i,x) -$ isod$_{C_i,x}(\kn^{ea}_{C/\P^2}, \ko_C) \le
0$, which  can be neglected in the right--hand side  of the $3d$--criterion (we
actually obtain
$-1$ if the node results from the intersection of two global components of
$C$).

\item Since, for  a node $\kn^{ea}_{C/S,x} = \kn^{es}_{C/S,x}$, $\tau_{C,x} =
\tau^{es}_{C,x}$ and since the isomorphism defect is a local invariant, we can
neglect nodes also in the $3d$--formulas for $es$.   For the $4d$--criterion,
however, nodes have to be counted with 1 (Lemma 4.7).
\end{enumerate}}
\end{remark}

\begin{sub}{\rm
 {\bf Examples}:
\begin{itemize}
\item[(b)] the irreducible curve $C \subset \P^2$ with affine equation
\vspace{-0.5cm}

\begin{eqnarray*}
f(x,y)& = & y^2 - 2x^2y + c_1 xy^2 +
c_2y^3+x^4-2c_1x^3y+c_3x^2y^2+c_4xy^3+c_5y^4+c_1x^5\\
& & \phantom{y^2} -(3c_2+2c_3)x^4y
- - -(2c_4+2)x^3y^2-(2c_1+2c_5)x^2y^3+c_6xy^4-(c_4+2)y^5\\
& & \phantom{y^2}+(2c_2+c_3)x^6+(c_4+2)x^5y+ (2c_1+c_5)
x^4y^2-c_6x^3y^3+(c_4+3)x^2y^4\\
& & \phantom{y^2}+ c_1xy^5-(3c_2+c_3+c_6)y^6
\end{eqnarray*}
where $c_1 := 16\alpha\beta^2 - 66\beta^2,\; c_2 := 3 \alpha\beta -
\frac{23}{2}\beta,\; c_3 := -(\beta + 7c_2),\; c_4 := -4\alpha + 13,\; c_5 :=
\beta^2-c_1,\; c_6 := 24\alpha\beta - 92\beta$ and $\alpha,\beta$ (complex)
solutions of $4\alpha^2 - 30\alpha + 55 = 0$ and $\beta^3 = \alpha^2 - 7\alpha
+ 12$ of degree 6 has exactly one singularity, which is of type $A_{19}$.
(The equation of this curve was found by H.\ Yoshihara (cf.\ \cite{Yos})).
Now $4d = 24 > 23 = \tau(C) + 4$ and, hence, each combination of
$A$--singularities given by an adjacent subdiagram of $A_{19}$ occurs on a
curve of degree 6 (this is a very simple proof of a well--known result which
was previously proved by using moduli theory of $K3$--surfaces).

\item[(c)] The curve $C \subset \P^2$ with homogeneous
equation $x^9 + zx^8 + z(xz^3+y^4)^2 = 0$ has exactly one  singularity  at
$(0:0:1)$ which is of type $A_{31}$.   Again we have  $4d=36 > 35 = \tau(C) +
4$, hence
$H^{ea}_{\P^2}$ is smooth at $C$.  $C$  is obtained by a small deformation of
Luengo's
example of a degree 9 curve ($x^9 + z(xz^3+y^4)^2=0$, having an
$A_{35}$--singularity) with non--smooth
$(H^{ea}_{\P^2}, C)$.   Of course, our criterion supports also this
non--smoothness since $4d = 36 < 39 = \tau(C) + 4$.

\item[(d)] $x^7 + y^7 + (x-y)^2 x^2y^2z = 0$ defines an irreducible curve $C
\subset \P^2$ which has 3 transverse cusps (not quasihomogeneous!) at
$(0:0:1)$ and no other singularities.   Since $4d = 28 > 24 + 4 - 1 \ge
\tau(C) + 4 -$ isod$_{(0:0:1)} (\widetilde{\kn}^{ea}_{C'/\P^2}, w_{C'})$ we see
that every local deformation of 3 transverse cusps can be realized by curves
of degree 7.
\end{itemize}
}
\end{sub}
\newpage


\end{document}